\documentclass[11pt, reqno]{amsart}
\usepackage{amsmath, amsthm, amscd, amsfonts, amssymb, graphicx, color}
\usepackage[bookmarksnumbered, colorlinks, plainpages]{hyperref}

\makeatletter \oddsidemargin.9375in \evensidemargin \oddsidemargin
\newtheorem{theorem}{Theorem}[section]
\newtheorem{lemma}[theorem]{Lemma}

\newtheorem{corollary}[theorem]{Corollary}
\theoremstyle{definition}
\newtheorem{definition}[theorem]{Definition}
\newtheorem{example}[theorem]{Example}

\newtheorem{remark}[theorem]{Remark}
\theoremstyle{approach}

\numberwithin{equation}{section}

\begin{document}
\setcounter{page}{1}

\title{$C^{\ast}-$algebra structure on vector valued-Banach algebras}

\author[Mitra Amiri and Ali {Rejali}$^*$]{Mitra Amiri and Ali {Rejali}$^*$}

\thanks{* Corresponding auther}


\begin{abstract}
Let $ \mathcal{A}$ be a commutative semisimple Banach algebra, $ X $ be a locally compact Hausdorff topological space and $ G $ be a locally compact topological group. In this paper, we investigate several properties of vector valued Banach algebras $ C_0(X , \mathcal{A})$, $ L^p(G,\mathcal A)$, $\ell^p (X , \mathcal{A})$ and $ \ell^{\infty}(X , \mathcal{A})$. We prove that these algebras are isomorphic with a $C^{\ast}-$algebra if and only if $ \mathcal{A} $ is so.
\end{abstract}

\maketitle

\section{Introduction and preliminaries}

Let $\mathcal{A}$ and $ \mathcal{B}$ be two Banach algebras. Then $\mathcal{A}$ is isomorphic with $\mathcal{B} $ and write
$ \mathcal{A} \cong \mathcal{B} $ if and only if there exists a map $ \varphi : \mathcal{A} \longrightarrow \mathcal{B}$ such that $\varphi$ is homomorphism, bijection and homeomorphism.

Let $ X $ be a locally compact Hausdorff topological space and $\mathcal{A}$ be a commutative Banach algebra such that the characteristic space $ \Delta( \mathcal{A})$ is non-empty. It has shown in \cite{S}, the conditions that $ C_0(X , \mathcal{A})$ is a BSE-algebra. In this paper, we investigate such conditions which $C_0(X,\mathcal{A})$ is homomorphic with a $C^{\ast}-$algebra. Here, we provide some preliminaries which will be used throughout the paper. We denote $ C_b(\Delta(\mathcal{A}))$, the Banach algebra consisting all continuous bounded complex valued functions on $\Delta(\mathcal{A})$ where $\Delta( \mathcal{A})$ is the character space of $\mathcal{A}$ containing all nonzero linear functional on $ \mathcal{A} $ which are multiplicative. We equipped $\Delta(\mathcal{A})$ with Gelfand topology. Clearly, $C_b(\Delta(\mathcal{A}))$ is a Banach algebra with pointwise product and supremum norm.

The function $\sigma\in C_b(\Delta(\mathcal{A}))$ is called BSE-function if there exists $\beta> 0$ such that for any complex number $c_1 , \cdots , c_n $ and $ \varphi_1 , \cdots , \varphi_n $ in $\Delta (\mathcal{A})$, the inequality
\begin{equation}\label{e1}
\big\vert \sum_{j=1}^{n} c_j \sigma ( \varphi_j ) \big\vert \leq  \beta \Vert \sum_{j=1}^{n} c_j \varphi_j \Vert_{\mathcal{A^{\ast}}}
\end{equation}
holds. The set of all BSE-functions is denoted by $C_{BSE}(\Delta(\mathcal{A}))$. Takahasi and Hatori in \cite{TH} showed that $C_{BSE}(\Delta(\mathcal{A}))$ is a commutative semisimple Banach algebra, equipped with the norm $ \Vert \cdot \Vert_{BSE} $ which is the infimum of all such $ \beta $ that satisfied \ref{e1}.
Note that,
$$
\Vert \sigma \Vert_{BSE} = sup \bigg\lbrace  \big\vert \sum_{j=1}^{n} c_j \sigma ( \varphi_j )\big\vert :\, \Vert \sum_{j=1}^{n} c_j \varphi_j \Vert_{\mathcal{A^{\ast}}} \leq 1 \bigg\rbrace.
$$
Consider the Gelfand mapping $\Gamma_{\mathcal{A}} : \mathcal{A} \rightarrow C_b (\Delta(\mathcal{A}))$ is defined by $a\mapsto\widehat{a},$ such that $\widehat{a}$ is the Gelfand
transform of $a$ and $\widehat{a}(\varphi)=\varphi(a)$ for $\varphi\in\Delta( \mathcal{A})$. Clearly, the Gelfand map $\Gamma_{\mathcal{A}}$ is continuous and homomorphism such that $\Vert\widehat{a}\Vert_{\infty}\leq\Vert a \Vert$ for all $a\in\mathcal{A}$. Furthermore, $\Gamma_{\mathcal{A}}$ is injective if and only if $\mathcal{A}$ is semisimple. Also $\Gamma_{\mathcal{A}}$ is bounded below if and only if $\Gamma_{\mathcal{A}}$ is renge closed, \cite{K}.

A commutative Banach algebra $\mathcal{A}$ is called without order if $a \mathcal{A}=\lbrace 0 \rbrace$ implies that $a = 0$, for $a \in \mathcal{A}$. It is easy to show that every commutative semisimple Banach algebra is without order. Let $\mathcal{A}$ be a commutative without order. Then a bounded linear operator on $\mathcal{A}$ is called a multiplier if it satisfies
$T(ab) = aT(b)$, for all $a,b \in \mathcal{A}$.
The set of all multipliers on $\mathcal{A}$ is denoted by $M(\mathcal{A})$ which is called the multiplier algebra of $\mathcal{A}$ and is a unital commutative
Banach algebra.
Following \cite[Theorem 1.2.2]{L}, for any $T\in M(\mathcal{A})$, there exists a unique bounded continuous function $\widehat{T}$ on $\Delta(\mathcal{A})$ such that
$$  \widehat{T(a)} ( \varphi ) = \widehat{T}( \varphi ) \hat{a} ( \varphi),$$ for
all $a\in {\mathcal A}$ and $\varphi\in \Delta({\mathcal A})$. Put
$$ \widehat{M( \mathcal{A} )} = \lbrace  \hat{T} :\, T \in  M( \mathcal{A} ) \rbrace $$
In \cite[Theorem 1.2.2]{L}, they show that $ \Vert \hat{T} \Vert_{\infty} \leq \Vert T \Vert $ for each $ T \in M( \mathcal{A} ) $. The Banach algebra $ \mathcal{A} $ is called \emph{BSE-algebra} if it satisfies the condition $$C_{BSE}(\Delta ( \mathcal{A}))= \widehat{M(\mathcal{A})}.$$
Let $ \mathcal{M} ( \mathcal{A} ) = \lbrace  \Theta \in C_b ( \Delta ( \mathcal{A})):\, \Theta \cdot \hat{\mathcal{A}} \subseteq \hat{\mathcal{A}} \rbrace $. Then
$  \widehat{M( \mathcal{A} )} \subseteq \mathcal{M} ( \mathcal{A} ).$
If $ \mathcal{A} $ is semisimple, then $ \widehat{M(\mathcal{A})}= \mathcal{M} ( \mathcal{A} ) $.

The BSE-algebra $ \mathcal{A} $ is called a type I-BSE-algebra, if
\begin{center}
$C_{BSE}(\Delta(\mathcal{A}))= C_b(\Delta(\mathcal{A}))$.
\end{center}
Takahashi and Hatori in \cite{TH} introduced and studied these algebras. They showed that a type I-BSE-algebra with a bounded approximate identity is isomorphic with a $C^{\ast}-$algebra and conversely.

It is known that
$$
 C_0 (X , \mathcal{A}) = C_0 (X) \check{\otimes} \mathcal{A} $$ as two Banach algebras which are isomorphic, where $ C_0 (X , \mathcal{A}) $ is the Banach algebra of all continuous functions $ f : X \longrightarrow \mathcal{A} $ which vanish at infinity, see \cite{K}.
Suppose, $ \Vert f \Vert_{\infty , \mathcal{A}} = \sup \big\lbrace  \Vert f(x) \Vert_{\mathcal{A}}:\, x \in X \big\rbrace$. Then $ C_0 (X , \mathcal{A} )$ is a commutative Banach algebra with pointwise product. In \cite{S}, the authors proved that the Banach algebra $ C_0(X ,\mathcal{A})$ is a BSE-algebra (type I-BSE-algebra) if and only if $ \mathcal{A}$ is BSE-algebra (type I-BSE-algebra), respectively.

 Let $ G$ be a locally compact Hausdorff topological group with Haar measure $ \lambda $, $\mathcal{A}$ be a separable commutative Banach algebra, $ 1 < p < \infty $ and $ L^p (G, \mathcal{A})$ be the Banach algebra of all Borel measurable function $ f : G \longrightarrow \mathcal{A} $ such that $ \Vert f \Vert_{p, \mathcal{A}} = \big( \int_G \Vert f(x) \Vert_{\mathcal{A}}^{p} d \lambda x\big)^{\frac{1}{p}} $ is finite. Then $ L^p (G, \mathcal{A})$ is a Banach algebra, under convolution product if and only if $ G $ is compact; see \cite{AAR}. In this paper, we show that $L^p(G,\mathcal{A})$ is a type I-BSE-algebra if and only if $ \mathcal{A} $ is a type I-BSE-algebra.

Dales in \cite{D} introduced BSE-norm algebra, if there exists $ M> 0 $ such that
\begin{center}
$ \Vert a \Vert_{\mathcal{A}} \leq M\Vert \hat{a} \Vert_{BSE} \;\;\;\;\;\;\;\;\;\;\;\; ( a \in \mathcal{A}).$
\end{center}
Clearly, $ \Vert \hat{a} \Vert_{BSE} \leq  \Vert a \Vert_{\mathcal{A}}$, for all $ a \in \mathcal{A} $. Thus, by using open mapping theorem the algebra $ \mathcal{A} $ is a BSE-norm algebra if and only if $ \mathcal{A} = \hat{\mathcal{A}} $ with equivalent norm as Banach algebras.

The algebra $C_0(X, \mathcal{A})$ [resp . $ L^p (G , \mathcal{A} ) ] $ is a BSE-norm algebra if and only if $ \mathcal{A} $ is a BSE-norm algebra, see \cite{S, AAR}.

Throughuot this paper we study the BSE-properties of the Banach algebra $\ell^p(X , \mathcal{A})$, where $1 \leq p \leq \infty $. Abtahi and Kamali in \cite{AK} show that $\ell^p (X , \mathcal{A})$ is BSE-algebra if and only if $\mathcal{A} $ is so for $1 \leq p \leq \infty$. In this paper, we investigate other properties of $\ell^p (X , \mathcal{A})$ as a type I-BSE-algebra and BSE-norm algebra. Furthermore, we show that $\ell^\infty(X , \mathcal{A})$ is a BSE-algebra (type I-BSE-algebra , BSE-norm algebra) if and only if $ \mathcal{A} $ is so, resepctively. Also, $ \ell^{\infty}(X ,\mathcal{A})$ is a $C^{\ast}-$algebra if and only if $ \mathcal{A} $ is a $ C^{\ast}-$algebra. In addition, $ M_b(G,w) $ is a type I-BSE algebra if and only if $ G $ is finite. Also the algebra $L^{1}(G, w)^{\ast \ast}$ is a type I-BSE algebra if and only if $ M_b(G, w)^{\ast \ast} $ is a type I-BSE algebra if and only if the group $ G $ is finite.

\section{Properties of $C_0(X, \mathcal{A})$}

In this section $ X $ is a locally compact Hausdorff topological space and $ \mathcal{A}$ is a commutative Banach algebra with non-empty character space. Let $g \in C_0 (X)$ and $a \in\mathcal{A}$. Then
$$
f = g \otimes a \in C_0(X ,\mathcal{A})
$$
such that
$$f(x) = g(x) a \;\;\;\;\;\;\;\;\;\;\;\; ( x \in X ).$$
Produces space by such functions are dense norm in $ C_0(X , \mathcal{A})$.
Suppose that $\mathcal{A}$ is a $\ast-$Banach algebra. Therefore, $C_0 (X , \mathcal{A})$ is a $\ast-$Banach algebra such that for all
$$ u = \sum_{j=1}^{n} g_j \otimes a_j \in C_0(X) \otimes \mathcal{A},$$
we define
$ u^{\ast} = \sum_{j=1}^{n} \overline{g_j} \otimes a_{j}^{\ast}.$

\begin{lemma}\label{l1}
Let $ C_0 (X , \mathcal{A} )$ be a $ C^{\ast}-$algebra. Then $\mathcal{A}$ is a $ C^{\ast}-$algebra.
\end{lemma}

\begin{proof}
Suppose that $a\in\mathcal{A}$ and $0 \neq g \in C_0(X)$. Put $f:= g \otimes a $. Thus
\begin{center}
$\Vert f \cdot f^{\ast}\Vert_{\infty , \mathcal{A}} = \Vert f \Vert_{\infty , \mathcal{A}}^{2}.$
\end{center}
So
\begin{align*}
\Vert f f^{\ast} \Vert_{\infty , \mathcal{A}} &= \sup \lbrace \Vert g(x) a \overline{g}(x) a^{\ast} \Vert_{\mathcal{A}} :\, x \in X \rbrace\\
&= \Vert a a^{\ast} \Vert_{\mathcal{A}} \,\sup \lbrace \vert (g\overline{g})(x) \vert :\, x \in X  \rbrace  \\
&= \Vert a a^{\ast} \Vert_{\mathcal{A}}\, \Vert g \overline{g} \Vert_{\infty} \\
&= \Vert a a^{\ast} \Vert_{\mathcal{A}} \, \Vert g \Vert_{\infty}^{2}.
\end{align*}
Moreover,
\begin{center}
$\Vert f \Vert_{\infty ,\mathcal{A}} = \Vert g \Vert_{\infty}\, \Vert a \Vert_{\mathcal{A}}.$
\end{center}
Consequently,
$$ \Vert a a^{\ast} \Vert_{\mathcal{A}} \, \Vert g \Vert_{\infty}^{2} =  \Vert a\Vert_{\mathcal{A}}^{2} \, \Vert g \Vert_{\infty}^{2}.
$$
Since $0 \neq g $, $ \Vert a a^{\ast} \Vert_{\mathcal{A}} = \Vert a\Vert_{\mathcal{A}}^{2}$. Therefore $\mathcal{A}$ is a $C^{\ast}-$algebra.
\end{proof}

\begin{lemma}\label{l2}
Let $\mathcal{A}$ be a commutative $C^{\ast}-$algebra. Then $C_0(X,\mathcal{A}) $ is a commutative $C^{\ast}-$algebra.
\end{lemma}

\begin{proof}
Consider $\mathcal{A} = C_0 (Y)$ where $Y= \Delta(\mathcal{A})$. In this case $ Y $ is a locally compact Hausdorff space and so $X \times Y $ is a locally compact Hausdorff space. Also
\begin{align*}
C_0(X ,\mathcal{A} ) &= C_0 (X) \check{\otimes} \mathcal{A} \\
&=C_0 (X) \check{\otimes} C_0 (Y) \\
&=C_0 (X \times Y),
\end{align*}
as two isometric Banach algebras.
Consequently, $C_0(X ,\mathcal{A})$ is a commutative $C^{\ast}-$algebra.
\end{proof}

By using lemmas \ref{l1} and \ref{l2} we have:

\begin{theorem}\label{T1}
Let $\mathcal{A}$ be a commutative $ C^{\ast}-$algebra and $ X $ be a locally compact Hausdorff space. Then $C_0 (X,\mathcal{A})$ is a $C^{\ast}-$algebra if and only if $\mathcal{A}$ is a $C^{\ast}-$algebra.
\end{theorem}

\begin{lemma}\label{l11}
Let $\mathcal{A}$ be a Banach algebra and $ X $ be a locally compact Hausdorff space. Then $C_0(X ,\mathcal{A})$ has a bounded approximate identity if and only if $\mathcal{A}$ has a bounded approximate identity.
\end{lemma}

\begin{proof}
Suppose that the net $(f_{\alpha} )_{\alpha}$ is a bounded approximate identity for $C_0 (X ,\mathcal{A})$ and $ x_0 \in X $ be constant. Put $e_{\alpha}:= f_{\alpha}(x_0)$. Therefore the net $(e_{\alpha})_{\alpha}$ is a bounded approximate identity for $\mathcal{A}$. In fact, if $g \in C_0(X)$ such that $ g(x_0)=1$, set $f_a := g \otimes a$ for $ a \in \mathcal{A}$. Hence, for all $\epsilon >0$, there exists $\alpha_0$ such that for all $\alpha \geq \alpha_{0}$ we have
$$\Vert f_{\alpha} \cdot f_{a} - f_{a} \Vert_{\infty , \mathcal{A}} < \epsilon.$$
Thus
\begin{align*}
\Vert e_{\alpha} a - a \Vert_{\mathcal{A}} &= \Vert f_{\alpha}(x_0)g(x_0)a-g( x_0)a \Vert_{\mathcal{A}} \\
&\leq \Vert f_{\alpha} \cdot f_{a} - f_{a} \Vert_{\infty , \mathcal{A}} < \epsilon.
\end{align*}
Therefore, $e_{\alpha} a \longrightarrow a$ for each $a \in \mathcal{A}$. Moreover,
\begin{center}
$\Vert e_{\alpha} \Vert_{\mathcal{A}}=\Vert f_{\alpha}(x_0) \Vert_{\mathcal{A}} \leq \Vert f_{\alpha} \Vert_{\infty,\mathcal{A}} < M$
\end{center}
for some $M> 0 $. So the net $(e_{\alpha})$ is a bounded approximate identity for $\mathcal{A}$.

Conversely, consider $ \epsilon > 0 $ and
$$ u = \sum_{j=1}^{n} g_j \otimes a_j \in C_0 (X , \mathcal{A}).$$
Put $\epsilon^{'}:= \frac{\epsilon}{(3 r + 1)n}$, where
$$  r = max \lbrace \Vert g_j \Vert_{\infty}\, \Vert a_j \Vert_{\mathcal{A}} :\, j = 1 , \cdots , n \rbrace.$$
Since $C_0(X)$ is a commutative $C^{\ast}-$algebra, it has a bounded approximate identity. Thus, there exists $g \in C_0(X)$ such that
\begin{center}
$ \Vert g_j - gg_j \Vert_{\infty} < \epsilon^{'} \;\;\;\;\;\;\;\;\;\;\;\; (j=1 , \cdots , n )$
\end{center}
Also
\begin{align*}
\Vert u - ( g \otimes a ) u \Vert_{\epsilon} &\leq \Vert u - ( g \otimes a ) u \Vert_{\pi}\\
&\leq \sum_{j=1}^{n} ( \Vert g_j -gg_j \Vert_{\infty}\, \Vert a_j \Vert_{\mathcal{A}} + \Vert g_j \Vert_{\infty}\, \Vert a_j -aa_j \Vert_{\mathcal{A}}
\\ &+ \Vert g_j - gg_j \Vert_{\infty}\, \Vert a_j - aa_j \Vert_{\mathcal{A}} \\
&< \epsilon
\end{align*}
such that $ \Vert \cdot \Vert_{\epsilon} $ and $ \Vert \cdot \Vert_{\pi} $ are injective and projective norms, respectively. Since $C_0(X) \otimes \mathcal{A} $ is dense in $C_0(X)\check{\otimes}\mathcal{A}$, according to \cite[Theorem 2.9.14]{D}, the commutative algebra $C_0(X, \mathcal{A})$ has a bounded approximate identity.
\end{proof}

\begin{theorem}\label{T2}
Let $\mathcal{A}$ be a commutative Banach algebra with a bounded approximate identity and $ X $ be a locally compact Hausdorff space. Then $C_0(X , \mathcal{A})$ is isomorphic with a $C^{\ast}-$algebra if and only if $ \mathcal{A}$ is isomorphic with a $C^{\ast}-$algebra.
\end{theorem}

\begin{proof}
Suppose that $C_0 (X , \mathcal{A})$ is isomorphic with a $C^{\ast}-$algebra. By \cite[Theorem 3]{TH}, $C_0 (X , \mathcal{A})$ is a type I-BSE algebra. Therefore according to \cite{S}, $\mathcal{A}$ is a type I-BSE algebra. Furthermore, according to assumption, $\mathcal{A}$ has a bounded approximate identity. Again, by \cite[Theorem 3]{TH}, $ \mathcal{A}$ is isomorphic with a $C^{\ast}-$algebra.

Conversely, let $\mathcal{A}$ be isomorphic with a $C^{\ast}-$algebra. Then by \cite[Lemma 2]{TH}, there exists $\beta > 0$ such that for any
$c_1 , \cdots , c_n \in \mathbb{T}$ and the same number of $\varphi_1 , \cdots , \varphi_n$ in $\Delta(\mathcal{A})$, there exists $a\in \mathcal{A}$ where
$\Vert a \Vert \leq \beta$ such that $\hat{a}(\varphi_i)= c_i$ $(i= 1, \cdots , n).$ Suppose that
\begin{center}
$ x_i \otimes \varphi_i\, \in\,\Delta(C_0 (X , \mathcal{A}))= X \otimes \Delta (\mathcal{A} )$
\end{center}
and put $K= \lbrace x_1 , \cdots , x_n \rbrace $. Then $ K\subseteq X$ is compact. Therefore by Urysohn's lemma there exists $g\in C_0(X)$ such that
$g \vert_{K} = 1$ and $\Vert g \Vert_{\infty} \leq 1 $. Now
put $ f:= g \otimes a $. In this case $ f \in C_0(X , \mathcal{A})$ and
$$
\Vert f \Vert_{\infty , \mathcal{A}} = \Vert g \Vert_{\infty} \cdot \Vert a \Vert_{\mathcal{A}} \leq \Vert a \Vert_{\mathcal{A}} \leq \beta.$$
Moreover,
$$
\hat{f}(x_i \otimes \varphi_{i})=\varphi_{i}(f(x_i))=\varphi_{i}(g(x_i)a)
= \varphi_{i}(a) = c_i \;\;\;\;\;\;\;\;\;\;\;\; (i = 1 , \cdots , n).$$
Hence, there exists $\beta > 0$ such that for any
$c_1,\cdots,c_n \in \mathbb{T}$ and the same number of $\psi_1 ,\cdots, \psi_n$ in $\Delta(C_0 (X ,\mathcal{A}))$, there exists $f\in C_0(X , \mathcal{A})$ where
$\Vert f \Vert \leq \beta$ such that $\hat{f}(\psi_i)= c_i$ $(i= 1,\cdots, n).$
Then, by \cite[Lemma 2]{TH}, $C_0(X, \mathcal{A})$ is isomorphic with a $ C^{\ast}-$algebra.
\end{proof}

\begin{theorem}
Let $\mathcal{A}$ be a commutative Banach algebra and $ X $ be a locally compact Hausdorff space. If  $\mathcal{A}$ has a bounded approximate identity, then $C_0(X, \mathcal A)$ is a type I-BSE algebra if and only if $ \mathcal A $ is a $C^{\ast}- $algebra.
\end{theorem}	

\begin{proof}
	Since $\mathcal{A}$
	 has a bounded approximate identity, by lemma \ref{l11}, $ C_0 (X, \mathcal A)$ has a bounded approximate identity. Moreover, if $C_0(X, \mathcal A)$ is a type I-BSE algebra. Then by \cite[Theorem 3]{TH},
	$C_0(X , \mathcal{A})$ is a $C^{\ast}-$algebra. Hence by lemma \ref{l1}, $ \mathcal{A}$ is a $C^{\ast}-$algebra.

	Conversely, let $ \mathcal A $ be a $C^{\ast}-$algebra. Then by lemma \ref{l2},
	$C_0(X , \mathcal A)$
	is a $ C^{\ast}-$algebra. Again by \cite[Theorem 3]{TH}, $C_0(X , \mathcal{A})$ is a type I-BSE algebra.
\end{proof}

\section{Properties of  $\ell^{p}(X, \mathcal{A})$}

let $\mathcal{A}$ be a commutative Banach algebra, $X$ be a non-empty set and
$1\leq p<\infty$. Let
$$\ell^{p}(X, \mathcal{A})=\big\{f:X\longrightarrow \mathcal A: \sum_{x\in X} \lVert f(x)\rVert^{p} <\infty \big\}.
$$
It is easily verified that $\ell^{p}(X, \mathcal{A})$ is a commutative Banach algebra, endowed with the norm
$$ \lVert f \rVert_{p}=\bigg (\sum_{x\in X} \lVert f(x)\rVert^{p}\bigg)  ^{\frac{1}{p}}\;\;\;\;\;\;\;\;(f\in \ell^{p}(X, \mathcal{A}))$$
and pointwise product. We prove that $\ell^{p}(X, \mathcal{A})$ is a type I-BSE-algebra, under certain conditions.

\begin{theorem}
	Let $ \mathcal A $ be a commutative Banach algebra and $ X $ be a non-empty set. Then $\ell^p (X, \mathcal A)$ is a type I-BSE algebra if and only if $ X $ is finite and $ \mathcal A $ is a type I-BSE algebra.
\end{theorem}

\begin{proof}
	First, suppose that $ X $ is finite and $ \mathcal A $ is a type I-BSE algebra. By \cite[Theorem 2.3]{AK},
	\begin{center}
		$\Delta \big(\ell^p (X, \mathcal A )\big)= X \times \Delta (\mathcal A ).$
	\end{center}
	Since $ \mathcal A $ is a BSE-algebra. Hence, there exists $ \beta > 0 $ such that for all $ c_1 , \cdots , c_n \in \mathbb T $ and $ \varphi_1 , \cdots , \varphi_n \in \Delta ( \mathcal A ) $ there exists $ a \in \mathcal A $ such that $\Vert a \Vert_{\mathcal A}\leq \beta $ and $ \hat a ( \varphi_i ) = c_i\, ( 1 \leq i \leq n ) $. Take $ f := 1 \otimes a $ where $f(x) = a $ for all $ x \in X $, then $ f \in \ell^p (X , \mathcal A)$ and for all $ \psi_i , \cdots , \psi_n \in \Delta (\ell^p (X , \mathcal A )),$ there exist $ \varphi_i \in \Delta (\mathcal A ) $ and $ x_i \in X $  such that $ \psi_i = x_i \otimes \varphi_i\, ( 1 \leq i \leq n ) $. Therefore
	\begin{align*}
	\hat f(\psi_i ) &= \psi_i (f)= x_i \otimes \varphi_i (f)\\
	&= \varphi_i (f(x_i)) = \varphi_i (a)\\
	&= \hat a (\varphi_i)\,\,\,\, (1 \leq i \leq n )
	\end{align*}
	and
	$$  \Vert f \Vert_{p , \mathcal A}^{p} = \sum_{x \in G} \Vert f(x) \Vert_{\mathcal A}^{p} = \Vert a \Vert_{\mathcal A}^{p} \cdot |X |= \Vert a \Vert_{\mathcal A}^{p}\;\;\;\;\;\;\;(x\in X).$$
	Thus
	\begin{center}
		$\Vert f \Vert_{p, \mathcal A} = \Vert a \Vert_{\mathcal A} \leq \beta.$
	\end{center}
	Consequently $ \ell^p (X , \mathcal A ) $ is a type I-BSE-algebra.
	
	Conversely, let $\ell^p (X, \mathcal A ) $ be a type I-BSE algebra. Then, $ \ell^p (X, \mathcal A ) $ is a BSE-algebra. Thus by \cite[Theorem 3.1]{AK}, $ \mathcal A $  is a BSE-algebra and $ X $ is finite. Now, we show that $ \mathcal A $ is a type I-BSE algebra.
	Since $\ell^p(X, \mathcal A)$ is a type I-BSE algebra, then by \cite[Lemma 2]{TH}, there exists $ \beta >0 $ such that for all $ c_1 , \cdots , c_n \in \mathbb T $ and $\psi_1 , \cdots , \psi_n \in \Delta(\ell^p (X , \mathcal A))$, there exists $ f \in \ell^p (X, \mathcal A)$ where $\Vert f \Vert_{p ,\mathcal A}\leq \beta$ and $\hat f(\psi_i)= c_i \; ( 1 \leq i \leq n )$. Moreover, for constant $ x_0 \in X $, put $\psi_i = x_0 \otimes \varphi_i $. Hence, $ \psi_i \in \Delta (\ell^p (X , \mathcal A))$ and
	\begin{align*}
	\hat f ( \psi_i )  = \psi_i (f) &= x_0 \otimes \varphi_i (f)\\
	&= \varphi_i (f(x_0 )) = \varphi_i (a)\\
	& = \hat a ( \varphi_i ) \,\,\,\, ( 1 \leq i \leq n )
	\end{align*}
	in which $a:=f(x_0)$ and
	\begin{center}
		$  \Vert a \Vert_{\mathcal A} = \Vert f(x_0 ) \Vert_{\mathcal A} \leq \Vert f \Vert_{1 , \mathcal A} \leq \Vert f \Vert_{p , \mathcal A} \leq \beta.$
	\end{center}
	Thus, there exists $ \beta > 0 $ such that for all $ c_1 , \cdots , c_n \in \mathbb T $ and $\varphi_1 ,\cdots,\varphi_n \in \Delta( \mathcal A)$, there exists $ a \in \mathcal A $ where $ \Vert a \Vert_{\mathcal A} \leq \beta $ and $ \hat a(\varphi_i)=c_i \;( 1 \leq i \leq n )$. This implies that $\mathcal A$ is a type I-BSE algebra and completes the proof.
\end{proof}

Now, we study the BSE-properties of Banach algebra $\ell^{\infty}(X, \mathcal{A})$.

Let
\begin{center}
	$\ell^{\infty}(X, \mathcal{A}) = \big\lbrace f : X \longrightarrow \mathcal{A} :\, f \, is \, bounded \big\rbrace$
\end{center}
and
\begin{center}
	$ \beta X = w^{\ast} - cl \lbrace \delta_x :\, x \in X \rbrace \subseteq \ell_{\infty}(X)^{\ast},$
\end{center}
where $ \beta X $ is the Stone-Cech compactification of $X$.
Then $ \beta X $ is a compact Hausdorff space and $\delta_x(g)= g(x)$ for $ g \in \ell^{\infty}(X)$.

\begin{lemma}
Let $\mathcal{A}$ be a Banach algebra. Then
\begin{center}
$ \ell^{\infty}(X,\mathcal{A})=C(\beta X )\check{\otimes }\mathcal{A}$
\end{center}
as two isometric Banach algebras.
\end{lemma}

\begin{proof}
Define
\begin{center}
$ \Phi:\ell^{\infty}(X ,\mathcal{A})\longrightarrow C(\beta X ,\mathcal{A}), \;\;f \rightarrow \overline{f}$
\end{center}
where $\overline{f} (z) = z ( f)$, for each $z \in \beta X$ and $f \in \ell^{\infty}(X ,\mathcal{A})$.
Clearly $\Phi$ is linear. Now, we show that $\Phi$ is injective. Suppose that $\overline{f_1} = \overline{f_2}$. Then
\begin{center}
$f_1 (x) = \overline{f_1}(\delta_x)=\overline{f_2}(\delta_x)= f_2 (x)\;\;\;\;\;\;\;\;\; (x \in X ).$
\end{center}
Therefore $f_1=f_2 $. Let $h \in C(\beta X ,\mathcal{A})$. Then $h \in C(\beta X) \check{\otimes}\mathcal{A}$.

First, let $h = \sum_{k=1}^{n} g_k \otimes a_k $ in which $g_k \in C(\beta X)$ and $ a_k \in \mathcal{A} $. We define $ f(x):= \sum_{k=1}^{n} g_k (\delta_x) a_k $ and we show that $\overline{f} = h$. Suppose that $ z \in \beta X $. So there exists a net $(x_{\alpha})_{\alpha}$ such that $ \delta_{x_{\alpha}} \overset{w^{\ast}}{\longrightarrow} z$. Therefore, $g_k( \delta_{x_{\alpha}} ) \longrightarrow g_k(z)$ for all $ k = 1 , \cdots , n $. Hence,
$$  h(z) = \sum_{k=1}^{n} g_k (z) a_k = \underset{\alpha}{\lim} \sum_{k=1}^{n} g_k (\delta_{x_{\alpha}} ) a_k $$
and
$$ \overline{f}(z) = z(f) = \underset{\alpha}{\lim} f( x_{\alpha}) = \underset{\alpha}{\lim} \sum_{k=1}^{n} g_k (\delta_{x_{\alpha}})a_k.$$
Consequently, $h(z) = \overline{f}(z)$ for all $z \in \beta X $.

Now, let $h_n \in C(\beta X)\otimes \mathcal{A} $ such that $ h_n  \overset{\Vert \cdot \Vert_{\epsilon}}{\longrightarrow} h$. We define $f_n (x) := h_n(\delta_x)$, for all $ x \in X$. Since
\begin{align*}
h(\delta_x ) &= \underset{n}{\lim}\, h_n(\delta_x )\\
&= \underset{n}{\lim} f_n(x)\\
&= f(x) = \overline{f}(\delta_x)  \;\;\;\;\;\;\;\;\;\;\;\; ( x \in X ),
\end{align*}
then $ \overline{f_n} \longrightarrow \overline{f}$ for some $ f $ in $\ell^{\infty}(X ,\mathcal{A})$.
Thus according to continuity $ \overline{f} $ and $ h $, we have
\begin{center}
$  h(z) = \underset{n}{\lim}\, h(\delta_{x_{\alpha}} ) = \underset{n}{\lim} \, f_n(\delta_{x_{\alpha}}) = \overline{f} (z)\;\;\;\;\;\;\;\;\;\;\;\; ( z\in \beta X ).$
\end{center}
This implies that $\overline{f} = h $.

We show that $\Phi $ is isometric. Suppose that $f \in \ell^{\infty}(X , \mathcal{A})$. Then
\begin{align*}
\Vert \overline{f} \Vert_{\infty , \mathcal{A}} &= \sup \lbrace \Vert \overline{f}(z) \Vert_{\mathcal{A}} :\, z \in \beta X \rbrace\\
&= \Vert \overline{f}(z_0) \Vert_{\mathcal{A}}
\end{align*}
for some $ z_0 $ in the compact set $\beta X$. Let $ \delta_{x_{\alpha}} \overset{w^{\ast}}{\longrightarrow} z_0$. Then
\begin{center}
$ \Vert \overline{f} \Vert_{\infty , \mathcal{A}} = \underset{\alpha}{\lim}\, \Vert f(x_{\alpha})\Vert_{\mathcal{A}} \leq \Vert f \Vert_{\infty , \mathcal{A}}.$
\end{center}
Moreover,
\begin{align*}
\Vert f \Vert_{\infty , \mathcal{A}}  &= \sup \lbrace \Vert f(x) \Vert_{\mathcal{A}} :\, x \in X \rbrace\\
&= \sup \lbrace \Vert \overline{f}(\delta_{x}) \Vert_{\mathcal{A}}:\, x \in X \rbrace \\
&\leq \Vert \overline{f} \Vert_{\infty ,\mathcal{A}}.
\end{align*}
Thus
 $$ \Vert \overline{f} \Vert_{\infty , \mathcal{A}} = \Vert f \Vert_{\infty , \mathcal{A}}.$$
Finally, we prove that $\Phi$ is homomorphism. Let $ f_1 , f_2 \in \ell^{\infty}(X , \mathcal{A})$. Then
\begin{center}
$  \overline{f_1 f_2}(z)= \underset{\alpha}{\lim} \, f_1 f_2(x_{\alpha}) = \underset{\alpha}{\lim} \, f_1 (x_{\alpha}) \underset{\alpha}{\lim} \, f_2 ( x_{\alpha} ) = \overline{f_1}(z) \overline{f_2}(z)$
\end{center}
such that $ \delta_{x_{\alpha}} \overset{w^{\ast}}{\longrightarrow} z$. This implies that $ \overline{f_1} \overline{f_2} = \overline{f_1 f_2} $. Consequently, $ \Phi $ is homomorphism.
\end{proof}

\begin{theorem}\label{th20}
Let $\mathcal{A}$ be a commutative Banach $\ast-$algebra with the non-empty character space $\Delta(\mathcal{A})$ and $ X $ be a discrete space. Then
 \begin{enumerate}
	\item[(i)] $ \ell^{\infty} (X , \mathcal{A})$ is a BSE-algebra if and only if $ \mathcal{A} $ is a BSE-algebra.
\item[(ii)] $ \ell^{\infty}(X ,\mathcal{A})$ is a type I-BSE algebra if and only if $ \mathcal{A} $ is a type I-BSE algebra.
\item[(iii)] $\ell^{\infty} (X , \mathcal{A})$ is a BSE-norm algebra if and only if $ \mathcal{A}$ is a BSE-norm algebra.
\item[(iv)] $\ell^{\infty}(X ,\mathcal{A})$ is isomorphic with a $C^{\ast}- $algebra if and only if $\mathcal{A}$ is isomorphic with a $C^{\ast}-$algebra.
\item[(v)] $ \ell^{\infty}(X ,\mathcal{A})$ is a $C^{\ast}-$algebra if and only if $ \mathcal{A} $ is a $ C^{\ast}-$algebra.
\end{enumerate}
\end{theorem}

\begin{proof}
By \cite{S}, (i), (ii) and iii are proved and by following Theorems \ref{T1} and \ref{T2}, (iv) and (v) is clear.
\end{proof}

\begin{lemma}
	Let $ \mathcal A $ be a unital commutative $ C^{\ast}-$algebra. Then $ \ell^{\infty}(X ,\mathcal A)$ is a type I-BSE algebra.
\end{lemma}

\begin{proof}
	Since $\mathcal A$ is a commutative $C^{\ast}-$algebra, $\mathcal A $ is a  type I-BSE-algebra. Therefore
	$$C_{BSE}(\Delta(\mathcal A ))=C_b(\Delta(\mathcal A)).$$
	Consequently, according to \cite{S},
	\begin{center}
		$C_{BSE}(\Delta(C(\beta X ,\mathcal A)))=C_b(\Delta ( C(X , \mathcal A)))$
	\end{center}
	that is
	\begin{center}
		$C_{BSE}(\Delta(\ell^{\infty} (X , \mathcal A)))= C_b(\Delta ( \ell^{\infty} (X , \mathcal A ))).$
	\end{center}
	Furthermore, $ \mathcal A$ is semisimple. Then, $\ell^{\infty}(X, \mathcal A)$ is semisimple and unital. Thus
	\begin{align*}
	M ( \ell^{\infty}(X , \mathcal A )) &= \mathcal M(\ell^{\infty}(X, \mathcal A )) = \ell^{\infty} (X, \mathcal A )\\
	&= C(\beta X , \mathcal A) = C(\beta X ) \check{\otimes} \mathcal A\\
	&= C( \beta X) \check{\otimes} C( \Delta ( \mathcal A )) = C ( \beta X \times \Delta (\mathcal A )) \\
	&= C_b (\Delta (\ell^{\infty}( X , \mathcal A )))= C_{BSE} ( \Delta (\ell^{\infty} (X, \mathcal A ))).
	\end{align*}
	Hence, $\ell^{\infty}(X, \mathcal A)$ is a BSE-algebra. Also, there exists $\beta > 0$ such that for all $c_1 , \cdots , c_n \in \mathbb T$ and $\varphi_1 , \cdots , \varphi_n \in \Delta(\mathcal A)$ there exists $a\in\mathcal A$ such that $\Vert a \Vert_{\mathcal A} \leq \beta $ and $\hat a(\varphi_i)=c_i \,\, ( 1 \leq i \leq n ) $. Take $\psi_i \in \Delta (\ell^{\infty} (X , \mathcal A ))$, then there exist $\varphi_i \in \Delta(\mathcal A)$ and $z_i \in \beta X $ such that $ \psi_i = z_i \otimes \varphi_i $. Put $ f = 1 \otimes a $ which $f(x)= a$ for all $ x \in X $. Therefore
	$$\Vert f \Vert_{\infty , \mathcal A} = \Vert a \Vert_{\mathcal A} < \beta$$
	and
	$$\hat f(\psi_i)= \psi_i (f) = z_i \otimes \varphi_i(f)=z_i( \varphi_i \circ f ).$$
	Moreover, $\delta_{x_{\alpha}}\overset{w^{\ast}}{\longrightarrow} z_i $, hence
	\begin{center}
		$ z_i ( \varphi_i \circ f ) = \underset{\alpha}{\lim}\, \varphi_i \circ f( x_{\alpha}) = \underset{\alpha}{\lim}\, \varphi_i (f( x_{\alpha} )) = \varphi_i (a) = c_i .$
	\end{center}
	This implies that $\hat f(\psi_i)= c_i $. Therefore by \cite[Lemma 2]{TH}, $ \ell^{\infty}(X, \mathcal A )$ is a type I-BSE algebra.
\end{proof}

\begin{remark}\label{re}
	Let $\ell^{\infty}(X , \mathcal A)$ be a type I-BSE algebra. Then $ \mathcal A$ is a type I-BSE algebra.
\end{remark}

\begin{proof}
	Since $\ell^{\infty}(X ,\mathcal A)=C(\beta X , \mathcal A)$, according to \cite{S}, $\mathcal A $ is a type I-BSE algebra.
\end{proof}

	\begin{theorem}
	Let $\mathcal A $ has a bounded approximate identity. Then $\ell^{\infty}(X , \mathcal A)$ is a type I-BSE algebra if and only if $\mathcal A $ is a $C^{\ast}-$algebra.
\end{theorem}

\begin{proof}
Suppose, $\ell^{\infty}(X , \mathcal A)$ is a type I-BSE algebra. By remark \ref{re},
$ \mathcal A$ is a type I-BSE algebra. Moreover, $\mathcal A $ has a bounded approximate identity. Thus by \cite[Theorem 3]{TH}, $\mathcal A $ is a $C^{\ast}-$algebra.

Conversely, since $\mathcal A $ is a $C^{\ast}-$algebra, again by \cite[Theorem 3]{TH},
$ \mathcal A$ is a type I-BSE algebra. Hence, by theorem \ref{th20}, $\ell^{\infty}(X , \mathcal A)$ is a type I-BSE algebra.
	\end{proof}	

\section{Properties of $L^1(G,w)$}

  Let $ G $ be a locally compact Hausdorff topological group. A weight on $ G $ is a Borel measurable function $w: G \longrightarrow \mathbb{R}^{+}$ such that
   $$ w(xy)\leq w(x)w(y)\,\,\,\,(x, y\in G).$$
  We define $\Omega:G \times G\rightarrow (0,1]$ by
  $\Omega(x,y)=\frac{w(xy)} {w(x)w(y)}$.

  A function $h:X \times Y\rightarrow \mathbb{C}$ is said to be $0$-cluster if
  $$\lim_{n}\lim_{m}h(x_{n},y_{m})=0=\lim_{m}\lim_{n}h(x_{n},y_{m})$$
for each two sequences $\{x_{n}\}\subseteq X$ and $\{y_{m}\}\subseteq Y$ of distinct points, provided the involved limits exist.

Define $w^*$ on $G$ by $ w^*(x)= w(x)w(x^{-1}),\,(x\in G)$. It can be simply verified that $w^*$ is also a weight on $G$. Moreover,
$w^*$ is bounded on $G$ if and only if $w$ is semi-multiplicative (that is, there exists $c>0$ such that $cw(xy)\leq w(x)w(y)$, for all $x, y\in G$). Then $\Omega$ can not be $0$-cluster when $w^*$ is bounded. The weighted algebra
$$L^1(G,w):=\big\{f:G\longrightarrow \mathbb C ,\; \text{f is measurable } \big\}$$
and
$$\lVert f\rVert_{1,w}:=\int_{G}\lvert f(x)\rvert w(x)d(x)< \infty.$$
In fact
 $$L^1(G,w)=\big\{f:fw\in L^1(G) \big\}$$ and for $E\in B(G)$,
$\mu w(E)=\int_{E}w d\mu$. If $w\geq 1$, then
$$M_b(G,w)=\big\{\mu\in M_b(G):\mu w\in M_b(G) \big\}$$ and
$\lVert \mu \rVert_{w}=\lVert \mu w \rVert_{1}.$

\begin{definition}
Let $ \mathcal A$ be a commutative Banach algebra and $M(\mathcal{A})$ be the multiplier algebra of $ \mathcal A$. Then
	$\widehat{M( \mathcal A)}$ is $\Delta-$weak closed, if $\widehat{T_{\alpha}} \overset{w^{\ast}}{\longrightarrow} \sigma $ for a net
	$\lbrace\widehat{T_{\alpha}} \rbrace_{\alpha} $ and a $\sigma \in C_b\big(\Delta(\mathcal A)\big)$, then there exists $T \in M(\mathcal A)$ such that $\sigma = \widehat{T} $ and conversely.
\end{definition}

\begin{theorem}
Let $ \mathcal A$ be a semisimple commutative Banach algebra. Then the following statements are equivalent:
\begin{enumerate}
	\item[(i)] $\mathcal A$ is a BSE-algebra.
	\item[(ii)] $\mathcal A$ has a bounded $\Delta-$weak approximate identity and $\widehat{M(\mathcal A)}$ is $\Delta-$weak closed.
\end{enumerate}
\end{theorem}

\begin{proof}
$(ii)\Rightarrow(i)$ Since $\mathcal A$ has a bounded $\Delta-$weak approximate identity, by \cite[Corollary 5]{TH}, $  \widehat{M ( \mathcal A) } \subseteq C_{BSE} ( \Delta ( \mathcal A)).$

Conversely, Let $\sigma \in C_{BSE}(\Delta(\mathcal A))$. By \cite[Theorem 4]{TH}, there exists a bounded net $(a_n)_{n}$ in $\mathcal A$ such that
$ \widehat{a_n} \overset{w^{\ast}}{\longrightarrow} \sigma$ on
 $\Delta(\mathcal A)$, that is, there exists $M > 0$ such that for each $n$, $\rVert a_n\lVert < M$ and for each $ \varphi \in \Delta ( \mathcal A)$,
$\varphi(a_n)\longrightarrow \sigma(\varphi)$. In other words, $ \widehat{ \mathcal A} \subseteq \widehat{M( \mathcal A)}$. Thus since $\widehat{a}_{n}\in \widehat{\mathcal A}$, then $ \widehat{a_n} \in \widehat{M( \mathcal A)}$ and $\widehat{a_n} \overset{ \sigma ( \widehat{M( \mathcal A)} , \Delta ( \mathcal A))}{\longrightarrow} \sigma $. Hence 	$$
\sigma \in \overline{\widehat{M( \mathcal A)}}^{w^{\ast}} = \widehat{M( \mathcal A)}$$
and this implies that
$ C_{BSE} ( \Delta ( \mathcal A)) \subseteq \widehat{M( \mathcal A)} $. This completes the proof and $\mathcal A $ is a BSE-algebra.

$(i)\Rightarrow(ii)$ Let $\mathcal A$ be a BSE-algebra. By \cite[Corollary 5]{TH},
 $\mathcal A$ has a bounded $\Delta-$weak approximate identity. Now we show that
$\widehat{M( \mathcal A)}$ is $\Delta-$weak closed. Suppose that
$\sigma \in\overline{\widehat{M( \mathcal A)}}^{w^{\ast}}$. By \cite[Theorem 4]{TH}, there exists a bounded net $\lbrace T_{\alpha} \rbrace $ in $ M( \mathcal A)$ such that $ \widehat{T_{\alpha}} (\varphi ) \longrightarrow \sigma (\varphi ) $, for each $\varphi \in \Delta (\mathcal A)$. Since $\lbrace T_{\alpha} \rbrace $ is bounded, then there exists $ k > 0 $ such that for each $\alpha$, $ \Vert T_{\alpha} \Vert \leq k$. Also $ \mathcal A$ is semisimple. Therefore there
exists $M > 0 $ such that $\Vert\widehat{T_{\alpha}}\Vert_{BSE} \leq M \Vert T_{\alpha}\Vert$. Now let $ c_1 , \cdots , c_n \in \mathbb{C} $ and   	
$ \varphi_1 , \cdots , \varphi_n \in \Delta (\mathcal A)$. We have
\begin{align*}
\vert   \sum_{i=1}^{n} c_i \sigma ( \varphi_{i} ) \vert &\leq \vert   \sum_{i=1}^{n} c_i \widehat{T_{\alpha}}  ( \varphi_{i} ) \vert + \vert   \sum_{i=1}^{n} c_i ( \widehat{T_{\alpha}}  ( \varphi_{i} )- \sigma ( \varphi_i )) \vert \\
&\leq Mk \Vert \sum_{i=1}^{n} c_i  \varphi_{i}  \Vert_{ {\mathcal A}^{\ast}} + \vert \sum_{i=1}^{n} c_i ( \widehat{T_{\alpha}} (\varphi_{i} )- \sigma ( \varphi_i )) \vert.
\end{align*}
Taking the limit with respect to $\alpha$, we obtain 	
	$$ \vert \sum_{i=1}^{n} c_i \sigma(\varphi_{i})\vert \leq Mk
	\Vert \sum_{i=1}^{n} c_i\varphi_{i}\Vert_{ {\mathcal A}^{\ast}}.$$
	This implies that $\sigma \in C_{BSE}(\Delta (\mathcal A))$ and
	 $\Vert\sigma \Vert_{BSE}\leq Mk$. In other words, $\widehat{M( \mathcal A)}= C_{BSE}( \Delta ( \mathcal A))$. Thus $\sigma \in \widehat{M(\mathcal A)}$, that is $\widehat{M( \mathcal A)}$ is $\Delta-$weak closed and the proof is complete.
	\end{proof}

\begin{example}
	Let $\mathcal A= L^{1}(G)$. Then $\mathcal A $ is a BSE-algebra. Then
	$$
		\widehat{M( \mathcal A)}= C_{BSE}(\Delta(\mathcal A)) = \widehat{M(G)} = \lbrace  \widehat{\mu}: \, \mu \in M(G)\rbrace.
	$$
	$\mathcal A$ has a bounded $\Delta-$weak approximate identity and  $\widehat{M( \mathcal A)}$ is $\Delta-$weak closed.
\end{example}

\begin{theorem}
	Let $ G $ be an abelian compact Hausdorff group, $ \mathcal A $ be a semisimple separable commutative and unital Banach algebra and $ 1 < p < \infty $. Then $L^p(G , \mathcal A)$ is a type I-BSE algebra if and only if $ G $ is finite and $ \mathcal A $ is a type I-BSE algebra.
\end{theorem}	

\begin{proof}
	First, suppose that $ \mathcal A $ is a type I-BSE algebra and $ G $ is a finite group. Therefore, $ \mathcal A $ is a BSE-algebra and $ G $ is finite. Hence, by \cite[Theorem 3.3]{AAR}, $L^p (G, \mathcal A)$ is a BSE-algebra. Since $ \mathcal A $ is a type I-BSE algebra, by \cite[Lemma 2]{TH}, there exists a positive number $ \beta $ such that for any finite number of $ c_1 , \cdots , c_n \in \mathbb T $ and the same number of $ \varphi_1 , \cdots , \varphi_n \in \Delta ( \mathcal A ) $, there exists $ a \in \mathcal A $ such that $ \Vert a \Vert \leq \beta $ and $ \hat{a} ( \varphi_i ) = c_i \,\, ( i = 1, \cdots, n )$. Now, take $ \chi_1 , \cdots , \chi_n \in \hat{G} $ and $ f = 1_{G} \otimes a $. We obtain
	\begin{align*}
	\hat{f}(\chi_i \otimes \varphi_i ) &= \chi_i(\varphi_i \circ f )\\
	&= \int_G \overline{\chi_i (x)} \varphi_i(1_G(x)a)dx\\
	&= \varphi_i (a) \int_G \overline{\chi_i(x)} dx \\
	&= \varphi_i (a) \chi_i ( 1_G ) \\
	&= \varphi_i (a)\\
	&= \hat{a}(\varphi_i) \;\;\;\;\;\;\;\;\;\;\;\; ( 1 \leq i \leq n ).
	\end{align*}
	Since $ L^p (G , \mathcal A ) $ is a BSE-algebra, again by \cite[Lemma 2]{TH}, $ L^p(G,\mathcal A)$ is a type I-BSE algebra.
	
	Conversely, suppose that $L^p(G,\mathcal A)$ is a type I-BSE algebra. Then $ L^p (G, \mathcal A ) $ is a  BSE-algebra, and so, $ G $ is finite and $ \mathcal A $ is a BSE-algebra. According to assumption,
	there exists $ \beta > 0 $ such that for any
	finite number of $ c_1 , \cdots , c_n \in \mathbb T $
	and the same number of
	$ \psi_1 , \cdots , \psi_n \in \Delta(L^p(G , \mathcal A))$, there exists $f \in L^p(G, \mathcal A)$ where $\Vert f \Vert_{p , \mathcal A} \leq \beta $ and $\hat{f}(\psi_j)= c_j \,\, ( j=1 \cdots n )$. If $ \varphi_1 , \cdots , \varphi_n \in \Delta ( \mathcal A)$ and $a=\sum_{x \in G}f(x)$, then $ a \in \mathcal A $ and
	\begin{align*}
	\hat{a} ( \varphi_j) = \varphi_j (a) &= \sum_{x \in G} \varphi_j ( f (x))\\
	&= \int_G \varphi_j \circ f(x) dx\\
	&= \int_G \overline{I(x)} \varphi_j \circ f (x) dx \\
	&= I \otimes \varphi_j(f)
	\end{align*}
	Where $ I(x) =1$, for each $ x \in G $.
	Take, $ \psi_j := I \otimes \varphi_j $. Hence,
	$$ \hat{a}( \varphi_j ) = \hat f ( \psi_j ) = c_j .$$
	Thus, $ \mathcal A $ is a type I-BSE algebra and
	$$  \Vert a \Vert_{\mathcal A} \leq \sum_{x \in G} \Vert f(x) \Vert_{\mathcal A} = \Vert f \Vert_{1, \mathcal A} \leq \Vert f \Vert_{p , \mathcal A} \leq \beta.$$
\end{proof}

\begin{theorem}
	Let $G$ be a compact group and $\mathcal A$ be a separable commutative Banach algebra. Then $L^{1}(G,\mathcal A)$ is type I-BSE algebra if and only if $G$ is finite and $\mathcal A$ is BSE-algebra.
\end{theorem}

\begin{proof}
	By \cite{AAR}, the proof is clear. 	
\end{proof}	

\begin{theorem}
	Let $ G $ be a locally compact Hausdorff topological abelian group and $w$ be a weight on $ G $. Then
 the following statements are equivalent
 \begin{enumerate}
 	\item[(i)] $L^1(G,w)$ is a type I-BSE algebra.
 	\item[(ii)] $L^1(G,w)$ is regular and amenable.
 	\item[(iii)] $ G $ is finite group.
 \end{enumerate}
\end{theorem}

	\begin{proof}
	$(iii) \Rightarrow (i)$ Suppose that $ G $ is a finite abelian group. Then by \cite [Theorem 4]{R}, $L^1(G,w)$ is a C$^*$-algebra and so by \cite[Theorem 3]{TH},
	$L^1(G,w)$ is a type I-BSE algebra.	
	
$(i) \Rightarrow (iii)$ Suppose that $L^1(G,w)$ is a type I-BSE algebra. Since
		$L^1(G,w)$ has a bounded approximate identity, by \cite[Theorem 3]{TH},	
		$L^1(G,w)$ is isomorphic with a commutative C$^*$-algebra. Thus, there exists a
		locally compact Hausdorff space	
		$X$	such that $L^{1}(G,w)\cong C_{0}(X)$ and by \cite [Theorem 4]{R}, $G$ is finite.
				
		By \cite[Theorem 4]{R}, (ii) is equivalent to (iii) and the proof is complete.
\end{proof}	

	\begin{lemma}	
		Let $ G $ be a topological group and the function $w$ be a weight on $ G $.  Then	
		\begin{enumerate}
			\item[(i)] The algebra $L^1(G,w)$ is a C$^*-$algebra if and only if $ G $ be a trivial group.
			\item[(ii)] $L^1(G,w)$ is isomorphic with a $C^{\ast}- $algebra if and only if $ G $ is a finite group.
		\end{enumerate}
	\end{lemma}	

	\begin{proof}
	(i) Suppose that $L^1(G,w)$ is a C$^*-$algebra. By \cite[Theorem 4]{R}, $w$ is
	 multiplicative. Hence, it can be assumed that $w=1$. Now, if there exists $x \in G$ where $x$ is non-identical, then
	$$f=i\delta_{x}+i\delta_{x^{-1}}+\delta_{y}+\delta_{y^{-1}}.$$
	Thus, $f^*=-i\delta_{x^{-1}}-i\delta_{x}+\delta_{y^{-1}}+\delta_{y}$. So,	
	$\Vert f.f^{*}\Vert_{w}=8$ and $\Vert f\Vert_{w}^{2}=16$. Consequently, $L^1(G,w)$ is not a C$^*-$algebra. This contradiction implies that $ G $ is a trivial group.	
	
	 (ii) Suppose that $L^1(G,w)$ is equal with a C$^*-$algebra. Then
it is regular. Also $L^1(G,w)$ has a bounded approximate identity	
	and it is weakly sequentially complete. Thus by \cite{D}, it has identity. This implies that $G$ is discrete. Moreover, the equality $\Vert\delta_{x}.\delta_{x}^*\Vert_{w}=\Vert\delta_{x}\Vert_{w}^{2}$ for every
	$x\in G$ implies that $w(x)=\Delta(x)^{\frac{1}{2}}$, for each $x\in G$ ($\Delta$ is the modular function of $G$), and this implies that $w$ is multiplicative. Now according to \cite[Theorem 2]{R},
	 $G$ is finite.

  Conversely, if $G$ is finite, then by \cite[Theorem 4]{R}, $L^1(G,w)$ is a C$^*-$algebra.
	\end{proof}	

		\begin{example}	
		If $G$ be discrete and $w\geq 1$, then $\ell^{1}(G,w)$ is norm dense in $\ell^{1}(G)$. Hence $\Delta(\ell^{1}(G,w))$
		is homeomorphic to $\Delta(\ell^{1}(G))$. Thus $\ell^{1}(G,w)$ is a BSE-algebra.
	\end{example}	

		\begin{remark}
		Let $ G $ be a locally compact group and $w\geq1$. Then $\Delta(M_{b}(G,w))$ is homeomorphic to $\Delta(M(G))$. Hence $M_{b}(G,w)$ is a BSE-algebra if and only if $G$ is discrete.
	\end{remark}
	
	\begin{theorem}
		Let $ G $ be a locally compact and $ w$ be Borel measurable weight function. Then $ M_b(G,w) $ is a type I-BSE algebra if and only if $ G $ is finite.
	\end{theorem}

		\begin{proof}
			Suppose that $ M_b(G,w) $ is a type I-BSE algebra, then $M_b(G,w)$ is unital and semisimple. Thus
			\begin{align*}
			M_b(G,w)  &= \mathcal M (M_b(G,w))\\
			&= C_{BSE}\big(\Delta ( M_b (G,w ))\big)\\
			&= C_b(\Delta( M_b(G,w))).
			\end{align*}
			Therefore $ M_b(G,w) $ is isomorphism with a commutative $C^{\ast}- $algebra with equivalent norms. Hence $ M_b(G,w)$ is 
	regular and amenable. Consequently $ L^{1}(G,w) $ is regular and amenable
			 and so, $G$ is finite.
			
			 Conversely, since
				$$ \mathbb{C}^{n} = \ell^{1}(G) = M_b (G,w) = M_b (G),$$
				$M_b(G,w)$ is a type I-BSE algebra.
		\end{proof}
	
	\begin{corollary}
		The algebra $L^{1}(G, w)^{\ast \ast}$ is a type I-BSE algebra if and only if $ M_b(G, w)^{\ast \ast} $ is a type I-BSE algebra if and only if the gorup $ G $ is finite.
	\end{corollary}
	
		\begin{proof}
			Let $ B $ be a Banach algebra such that $ B^{\ast \ast} $ is regular and amenable. Then $ B $ is regular and amenable. Therefore if $ L^{1}(G,w)^{\ast \ast} $ or $ M_b (G ,w)^{\ast \ast} $ is a type I-BSE-algebra, then $ L^{1}(G,w) $ or $M_b(G,w)$ is regular and amenable, respectively. Consequently, $ G $ is finite.
		\end{proof}

\vspace{9mm}

{\footnotesize \noindent
	M. Amiri\\
	Department of Pure Mathematis\\
	Faculty of Mathematics and Statistics\\
	University of Isfahan\\
	Isfahan, 81746-73441\\
	 Iran\\
	mitra75amiri@gmail.com\\
	m.amiri@sci.ui.ac.ir\\
	
	\noindent
	A. Rejali\\
	Department of Pure Mathematis\\
	Faculty of Mathematics and Statistics\\
	University of Isfahan\\
	Isfahan, 81746-73441\\
	Iran\\
	rejali@sci.ui.ac.ir\\
	a.rejali20201399@gmail.com

\end{document}